\documentclass[11pt,a4paper,reqno]{amsart}

\usepackage{amscd,amssymb,amsthm}
\usepackage{graphicx}

\usepackage{enumitem}
\usepackage[pdftex]{hyperref}
\usepackage[pdftex]{bookmark}
\hypersetup{breaklinks=true,
            colorlinks=false,
            citecolor=green,
            urlcolor=blue,
            linkcolor=blue,
            pdftitle={A Generalization of Riemann Sums},
            pdfauthor={Omran Kouba}}

\usepackage{enumerate}
\newtheorem{theorem}{Theorem}

\newcommand{\ds}{\displaystyle}
\newcommand{\bg}{\bigskip\goodbreak}
\newcommand{\vf}{\varphi}
\newcommand{\abs}[1]{\left\vert #1\right\vert }
\newcommand{\N}[1]{\muskip0=-2mu{\left|\mkern\muskip0\left|#1\right|\mkern\muskip0\right|}}

\title[A Generalization of Riemann Sums ]
{A Generalization of Riemann Sums}
\author{Omran Kouba}
\address{Department of Mathematics \\
Higher Institute for Applied Sciences and Technology\\
P.O. Box 31983, Damascus, Syria.}
\email{\href{mailto:omran_kouba@hiast.edu.sy}{omran\_kouba@hiast.edu.sy}}
\thanks{\textsc{\Small{Mathematical Reflections  \textbf{1} (2010)}}}
\begin{document}
\begin{abstract}
We generalize the property  that  Riemann sums of a continuous function corresponding to
equidistant subdivision of an interval converge to the integral of that function, and we give some
applications of this generalization. 
\end{abstract}
\bg

\maketitle
In \cite{Lu} Cezar Lupu asked the following question : \bg
\noindent Prove that
\begin{equation}\label{E1}
\lim_{n\to\infty}\sum_{k=1}^n\frac{\arctan \frac{k}{n}}{n+k}\cdot\frac{\varphi(k)}{k}=\frac{3\log2}{4\pi},
\end{equation}
where $\vf$ denotes the Euler totient function. In this short note we will prove the following theorem, that will, in particuler, answer this question.
\bg
\begin{theorem}\label{th1}
 Let $\alpha$ be a positive real number, and let $(a_n)_{n\geq1}$ be a sequence of positive real numbers such that
\[
\lim_{n\to\infty}\frac{1}{n^\alpha}\sum_{k=1}^na_k=L.
\]
Then, for every continuous function $f$ on the interval $[0,1]$, we have
\[\lim_{n\to\infty}\frac{1}{n^\alpha}\sum_{k=1}^nf\left(\frac{k}{n}\right)\,a_k=L\int_0^1\alpha x^{\alpha-1}f(x)\, dx.\]
\end{theorem}
\bg
\begin{proof}
We use the following two facts:\bg

\begin{description}
\item[{\it fact 1}] for $\beta>0$ we have $\ds \lim_{n\to\infty}\frac{1}{n^{\beta+1}}\sum_{k=1}^nk^{\beta}=\frac{1}{\beta+1}$.
\item[{\it fact 2}] if $(\lambda_n)_{n\geq 1}$ is a real sequence that converges to $0$, and $\beta>0$ then
$\ds \lim_{n\to\infty}\frac{1}{n^{\beta+1}}\sum_{k=1}^nk^{\beta}\lambda_k=0$.
\end{description}

Indeed, {\it fact 1.} is just the statement that the Riemann sums of the function $x\mapsto x^{\beta}$ corresponding to equidistant
subdivision of the interval $[0,1]$, converge to $\int_0^1x^{\beta}\,dx$.\bg
The proof of  {\it fact 2.} is an argument ``\`a la Ces\'aro''.
Since $(\lambda_n)_{n\geq 1}$ converges to $0$ it must be bounbed, and if we define
$\Lambda_n=\sup_{k\geq n}\abs{\lambda_k}$, then $\lim_{n\to\infty}\Lambda_n=0$. But, for $1<m<n$ we have
\begin{align*}
\abs{
\frac{1}{n^{\beta+1}}\sum_{k=1}^{n}k^{\beta}\lambda_k}
&\leq\frac{1}{n^{\beta+1}}\sum_{k=1}^{m}k^{\beta}\abs{\lambda_k}+
\frac{1}{n^{\beta+1}}\sum_{k=m+1}^{n}k^{\beta}\abs{\lambda_k}\\
&\leq\frac{m^{\beta+1}}{n^{\beta+1}}\Lambda_1+\Lambda_m.
\end{align*}
Let  $\epsilon$ be an arbitrary positive number. There is an $m_\epsilon>0$ such that $\Lambda_{m_\epsilon}<\epsilon/2$, then
we can find $n_\epsilon>m_\epsilon$ such that for every $ n>n_\epsilon$ we have $m_{\epsilon}^{\beta+1} \Lambda_1/n^{\beta+1}<\epsilon/2$.
Thus, 
\[ n>n_\epsilon\quad\Longrightarrow\quad
\abs{
\frac{1}{n^{\beta+1}}\sum_{k=1}^{n}k^{\beta}\lambda_k}<\epsilon.\] 
This ends the proof of {\it fact 2.}
\bg
Now, we come to the proof of our Theorem.
We start by proving the following property by induction on $p$.
\begin{equation}\label{E2}
\forall\,p\geq0,\qquad\lim_{n\to\infty}\frac{1}{n^{\alpha+p}}\sum_{k=1}^nk^pa_k=\frac{\alpha}{\alpha+p}L
\end{equation}

The base property ($p=0$) is just the hypothesis. So, let us assume that this is true for a given $p$, and
let 
\[\lambda_n= \frac{1}{n^{\alpha+p}}\sum_{k=1}^nk^pa_k-\frac{\alpha L}{\alpha+p},\]
(with the convention $\lambda_0=0$,) so that $\lim_{n\to\infty}\lambda_n=0$. Clearly, we have
\[
k^pa_k=k^{\alpha+p}\lambda_k-(k-1)^{\alpha+p}\lambda_{k-1}+\frac{\alpha L}{\alpha+p}\left(k^{\alpha+p}-(k-1)^{\alpha+p}\right),
\]
hence,
\begin{align*}
k^{p+1}a_k& =k^{\alpha+p+1}\lambda_k-k(k-1)^{\alpha+p}\lambda_{k-1}+\frac{\alpha L}{\alpha+p}\left(k^{\alpha+p+1}-k(k-1)^{\alpha+p}\right),\\
& =k^{\alpha+p+1}\lambda_k-(k-1)^{\alpha+p+1}\lambda_{k-1}+\frac{\alpha L}{\alpha+p}\left(k^{\alpha+p+1}-(k-1)^{\alpha+p+1}\right)\\
&\qquad\qquad\qquad\qquad-(k-1)^{\alpha+p}\lambda_{k-1}
-\frac{\alpha L}{\alpha+p}(k-1)^{\alpha+p}
\end{align*}
It follows that
\begin{equation*}
\frac{1}{n^{\alpha+p+1}}\sum_{k=1}^nk^{p+1}a_k=\lambda_n
-\frac{1}{n^{\alpha+p+1}}\sum_{k=1}^{n-1}k^{\alpha+p}\lambda_{k}+\frac{\alpha L}{\alpha+p}\left(1
-\frac{1}{n^{\alpha+p+1}}\sum_{k=1}^{n-1}k^{\alpha+p}\right).
\end{equation*}
Using {\it fact 1} and  {\it fact 2} we conclude immediately that
\begin{equation*}
\lim_{n\to\infty}\frac{1}{n^{\alpha+p+1}}\sum_{k=1}^nk^{p+1}a_k=\frac{\alpha L}{\alpha+p}\left(1-\frac{1}{\alpha+p+1}\right)=\frac{\alpha L}{\alpha+p+1}.
\end{equation*}
This ends the proof of \eqref{E2}.\bg
For a continuous function $f$ on the interval $[0,1]$, we define
$$I_n(f)=\frac{1}{n^\alpha}\sum_{k=1}^nf\left(\frac{k}{n}\right)a_k,\quad
\hbox{and}\quad J(f)=L\int_0^1\alpha x^{\alpha-1}f(x)\,dx.$$
Now, if $X^p$ denotes  the function $t\mapsto t^p$, then $(2)$ is equivalent to the fact
that
$\lim_{n\to\infty}I_n(X^p)=J(X^p)$, for every nonnegative integer $p$.
Using linearity, we conclude that $\lim_{n\to\infty}I_n(P)=J(P)$ for every polynomial function $P$.\bg
On the other hand, If $M=\sup_{n\geq1}\frac{1}{n^\alpha}\sum_{k=1}^na_k$ then $L\leq M$, and we observe that for every continuous functions $f$ and $g$ on $[0,1]$, and every positive integer $n$, we have
\begin{align*}
\abs{I_n(f)-I_n(g)}&\leq M\sup_{[0,1]}\abs{f-g}\\
\noalign{\text{and}}\\
\abs{J(f)-J(g)}&\leq M \sup_{[0,1]}\abs{f-g}.
\end{align*}

So, consider a continuous function $f$ on $[0,1]$. Let $\epsilon$ be an arbitrary positive number. Using Weierstrass Theorem there is a polynomial 
$P_\epsilon$ such that \[\N{f-P_\epsilon}_\infty=\sup_{x\in[0,1]}\abs{f(x)-P_\epsilon(x)}<\frac{\epsilon}{3M}\] 
Moreover, since
$\lim_{n\to\infty}I_n(P_\epsilon)=J(P_\epsilon)$, there exists a positive integer $n_\epsilon$ such that
$\abs{I_n(P_\epsilon)-J(P_\epsilon)}<\frac{\epsilon}{3}$ for every $n>n_\epsilon$. Therefore, for $n>n_\epsilon$, we have
$$
\abs{I_n(f)-J(f)}\leq \abs{I_n(f)-I_n(P_\epsilon)}+\abs{I_n(P_\epsilon)-J(P_\epsilon)}
+\abs{J(P_\epsilon)-J(f)}<\epsilon.$$
This ends the proof of Theorem \ref{th1}.
\end{proof}
\bg
\noindent{\bf Applications.}\bg
\begin{itemize}
\item  It is known that the Euler totient function $\vf$ has very erratic behavior, but on the mean we have the following beautiful result, see \cite[\S 18.5]{Har},
\begin{equation}\label{E3}
\lim_{n\to\infty}\frac{1}{n^2}\sum_{k=1}^n\vf(k)=\frac{3}{\pi^2}.
\end{equation}
Using Theorem \ref{th1}.  we  conclude that, for every continuous function $f$ on $[0,1]$ we have
\begin{equation}\label{E4}
\lim_{n\to\infty}\frac{1}{n^2}\sum_{k=1}^n f\left(\frac{k}{n}\right)\vf(k)=\frac{6}{\pi^2}\int_0^1xf(x)\,dx.
\end{equation}
Choosing $f(x)=\frac{\arctan(x)}{x(1+x)}$ we conclude that
\begin{equation}\label{E5}
\lim_{n\to\infty}\sum_{k=1}^n \frac{\arctan(k/n)}{k(n+k)}\vf(k)=\frac{6}{\pi^2}\int_0^1\frac{\arctan(x)}{1+x}\,dx.
\end{equation}
Thus, we only need to evaluate the integral $I=\int_0^1\frac{\arctan(x)}{1+x}\,dx$. The ``easy'' way to do this,
is to make the change of variables $x\leftarrow \frac{1-t}{1+t}$ to obtain
\begin{align*}
I&=\int_0^1\arctan\left(\frac{1-t}{1+t}\right)\frac{dt}{1+t}=\int_0^1\left(\frac{\pi}{4}-\arctan(t)\right)\frac{dt}{1+t}\cr
&=\frac{\pi}{4}\int_0^1\frac{dt}{1+t}-I
\end{align*}
Hence, $I=\frac{\pi}{8}\log 2$. Replacing back in \eqref{E5} we obtain \eqref{E1}.
\bg
\item  In the same line of ideas, If $\sigma(n)$ denotes the sum of divisors of $n$, then we have also
the following result, see \cite[\S 18.3]{Har},
$$\lim_{n\to\infty}\frac{1}{n^2}\sum_{k=1}^n\sigma(k)=\frac{\pi^2}{12}.$$
Using Theorem 1. we  conclude that, for every continuous function $f$ on $[0,1]$ we have
$$
\lim_{n\to\infty}\frac{1}{n^2}\sum_{k=1}^n f\left(\frac{k}{n}\right)\sigma(k)=\frac{\pi^2}{6}\int_0^1xf(x)\,dx.
$$
Choosing for instance $f(x)=\frac{1}{1+ax^2}$ we conclude that
$$\lim_{n\to\infty}\sum_{k=1}^n \frac{\sigma(k)}{n^2+ak^2}=\frac{\pi^2}{12 a}\log(1+a).$$
\bg
\item Starting from
$$ \lim_{n\to\infty}\frac{1}{n}\sum_{k=1}^n\frac{\vf(k)}{k}=\frac{6}{\pi^2},$$
which can be provrd in the same way as (3), we conclude
that, for every $\alpha\geq0$ we have
\begin{equation}\label{E6} \lim_{n\to\infty}\frac{1}{n^{\alpha+1}}\sum_{k=1}^nk^{\alpha-1}\vf(k)=\frac{6}{\pi^2(1+\alpha)}
\end{equation}
Also,
\begin{align*}
\lim_{n\to\infty}\frac{1}{n^{\alpha+1}}\sum_{k=1}^nk^{\alpha-1}\log(k/n)\vf(k)
&=\frac{6}{\pi^2}\int_0^1x^\alpha\log(x)\, dx\\
&=-\frac{6}{\pi^2(\alpha+1)^2}.
\end{align*}
Hence, using \eqref{E6} we obtain for $\alpha\geq0$ the following result :
\[
\frac{1}{n^{\alpha+1}}\sum_{k=1}^nk^{\alpha-1}\log k\,\vf(k)
=\frac{6\big((1+\alpha)\log n
-1\big)}{\pi^2(1+\alpha)^2}+o(1)
\]
\end{itemize}


\bigskip

\begin{thebibliography}{9}
\setlength{\itemsep}{5pt}

\bibitem{Lu}
C. Lupu, 
Problem U131, {\em Mathematical Reflections.\/}
~{\bf (4)}
(2009).
\bibitem{Har}
G. H. Hardy  \textsc{and}  E. M.Wright, 
An Introduction to the Theory of Numbers (5th ed.), {\em Oxford University Press.\/} 
(1980).

\end{thebibliography}
\end{document}